\documentclass[11pt]{article}
\usepackage{a4wide,amsfonts,amssymb,amsmath,amsthm, amssymb}
\def\proof{{\bf Proof\quad}}
\def\beginpf{\proof}
\def\qed{\hfill\rule{2.2mm}{2.2mm}\vspace{1ex}}
\def\endpf{\qed}

\newtheorem{theorem}{Theorem}[section]

\newtheorem{example}[theorem]{Example}
\newtheorem{lemma}[theorem]{Lemma}

\newtheorem{remark}[theorem]{Remark}

\newcommand{\f}{\mathbf{f}}
\newcommand{\dist}{\mathop{\rm dist}\nolimits}

\def\epsilon{\varepsilon}

\def\B{\mathbf{B}}

\def\CC{\mathbb C}

\def\J{\mathcal J}

\def\NN{\mathbb N}
\def\RR{\mathbb R}

\def\LL{\mathcal L}

\def\BB{\mathcal{B}}

\newcommand{\carl}{\mathop{\rm Carl}\nolimits}

\newcommand{\re}{\mathop{\rm Re}\nolimits}
\newcommand{\im}{\mathop{\rm Im}\nolimits}
\newcommand{\diag}{\mathop{\rm diag}\nolimits}

\def\text{\mbox}

\def\beq{\begin{equation}}
\def\eeq{\end{equation}}

\title{Weighted multiple interpolation and the control of perturbed semigroup systems}
\author{{Birgit Jacob\thanks{Fachbereich C - Mathematik und Naturwissenschaften, 
Arbeitsgruppe Funktionalanalysis,
Bergische Universit\"at Wuppertal,
Gau{\ss}stra{\ss}e 20,
D-42097 Wuppertal, Germany. \tt jacob@math.uni-wuppertal.de}
\qquad Jonathan R.~Partington\thanks{School of Mathematics,
University of Leeds,
Leeds LS2 9JT, U.K. \tt
J.R.Partington@leeds.ac.uk}
\qquad Sandra Pott\thanks{Faculty of Science,
Centre for Mathematical Sciences,
Lund University,
22100 Lund, Sweden. \tt sandra.pott@math.lu.se}}}

\begin{document}
\maketitle

\begin{abstract}
In this paper the controllabillity and admissibility of perturbed semigroup systems are studied, using tools from the theory of interpolation and Carleson measures. In addition, there are new results on the perturbation of Carleson measures and on the weighted interpolation of functions and their derivatives in Hardy spaces, which are of interest in their own right.
\end{abstract}


{\bf 2000 Subject Classification.} 30E05, 30D55, 30H10, 47A57, 93B05.

\section{Introduction}

Our general situation is that we have an infinite-dimensional linear system governed by the state equations:
\beq\label{eq:state}
{dx(t) \over dt}=Ax(t)+Bu(t), \qquad \hbox{with} \quad x(0)=x_0, \quad \hbox{say}.
\eeq
Here $u$ is the   input, and $x$ the state. The operator $A$ is assumed to be the
infinitesimal generator of a $C_0$ semigroup $(T(t))_{t \ge 0}$ on a Hilbert space $H$.
In general $B$  (the control  
operator) is unbounded, and in the case of a one-dimensional input space it
corresponds to a sequence $(b_n)$ as will be explained in more detail in Section \ref{sec:new2}.\\

This paper addresses two issues. First, in Section \ref{sec:new2} we consider controllability (or, dually, observability) of a system involving a perturbation of the generator $A$. The results here
are derived in the case of diagonal systems with scalar inputs and are obtained by means of an estimate involving the
perturbation of atomic Carleson measures.\\

 In Section \ref{sec:new3} we study weighted multiple interpolation: this  technique allows us to consider the controllability of systems with finitely-many
Jordan blocks (and a Riesz basis of generalized eigenvectors), extending the main results of
\cite{JP06}.
More sophisticated higher-order interpolation results are derived in Section  \ref{sec:new4}, which extend some classical results of Vasyunin \cite{nik}. These can be applied to analyse the controllability of systems with infinitely-many nontrivial Jordan blocks, provided that the sizes of the blocks are uniformly bounded.\\

The second issue addressed is admissibility, which is treated in Section \ref{sec:new5}. Here perturbation results for contraction semigroups
are obtained by linking admissibility with a resolvent condition (a particular case of the Weiss conjecture which was proved in \cite{JP01}).\\

A general reference for the material on controllability and admissibility is the book \cite{TW}, although we shall  define the necessary concepts
as they are needed.
We mention also some related work: in the paper \cite{GX} the exact controllability  of $C_0$ groups with multiple eigenvalues is analysed; 
some related work on perturbation
can be found in \cite{HI}, where it is shown  that admissibility of observation operators is invariant under Miyadera--Voigt perturbations of the operator $A$; in \cite{Hadd05} results are given on controllability and admissibility of systems with Desch--Schappacher perturbations.


\section{Controllability, and perturbations of Riesz bases}\label{sec:new2}

Our general framework here is that we consider an exponentially stable system given by (\ref{eq:state}), where the
operator $A$ has a Riesz basis $(\phi_n)$ of eigenvectors with eigenvalues $(\lambda_n)$ lying in the left half-plane $\CC_-$. We are concerned with
the persistence of controllability properties under perturbations of $A$, in particular,
those in which the perturbed operator also has a Riesz basis of eigenvectors. Later we shall also allow for the presence of nontrivial Jordan blocks.

The control operator $B: \CC \to \widetilde{H}$ given by $Bz=z\sum_n b_n \phi_n $ may be considered as an operator mapping into a larger space $\widetilde H$, the completion of the span of $(\phi_n)$ under a smaller norm, for example
\[
\left \| \sum_n a_n \phi_n\right  \|^2_1 = \sum_n \frac{|a_n|^2}{n^2(1+|b_n|^2)}. 
\]
The controllability operator $\BB_\infty: L^2(0,\infty) \to \widetilde H$ is defined by
\[
\BB_\infty u = \int_0^\infty T(t)Bu(s) \, ds,
\]
where we interpret the  semigroup $(T(t))_{t \ge 0}$ as acting on $\widetilde H$ in the obvious way. The system is said to be
{\em exactly controllable} if $R(\BB_\infty) \supset H$, and {\em null-controllable in time $\tau>0$} if $R(\BB_\infty) \supset R(T(\tau))$. Here $R(.)$ denotes the range of an operator.\\

There are various relevant perturbation results in the literature, the most recent being in
\cite{wyss}. For example, there is \cite[Thm. 4.15a]{kato} on perturbations of self-adjoint operators
with compact resolvents (but the self-adjoint condition is rather strong), and
\cite[Thm. XIX.2.7]{DS3} on perturbations of {\em discrete spectral operators\/} (these have compact resolvent,
a spectral resolution, and only finitely many repeated eigenvalues). One of the earliest perturbation results
is the following from \cite{DS3}.

\begin{theorem}[Dunford and Schwartz]
Let $T$ be Hilbert space discrete spectral operator with eigenvalues $(\lambda_n)$ and such that
all but finitely many of the spectral projections $E(\lambda_n,T)$ have one-dimensional range. Let
$d_n=\dist(\lambda_n, \sigma(T)\setminus \{\lambda_n\})$, and take $\lambda_0 \in \CC \setminus \sigma(T)$.
Take $0 \le \nu < 1$.
If $P$ is an operator such that $P(T-\lambda_0 I)^{-\nu}$ is bounded, and
\[
\sum_n \frac{(|\lambda_n| + d_n)^{2\nu}}{d_n^2} < \infty,
\]
then $T+P$ is a discrete spectral operator and $\sigma(T+P)$ consists of
a sequence $(\mu_n)$ with 
\[
|\lambda_n - \mu_n| \le C |\lambda_n|^\nu,
\]
and again all but finitely many of the spectral projections $E(\mu_n,T+P)$ have one-dimensional range.
\end{theorem}

For example, if $\lambda_n$ is of order $n^2$, and $d_n$ of order $n$, then
we require $4\nu - 2 < -1$ or $\nu < 1/4$.\\


The controllability of a diagonal system with scalar input was linked with a certain Carleson measure property
in \cite{JP06}; a key fact to recall is that the system is exactly controllable if and only if
the interpolation problem  $b_n F(-\lambda_n)= a_n$ for $n=1,2,\ldots$ can be solved for $F$ 
in the Hardy space $H^2(\CC_+)$ on the right half-plane
for all sequences $(a_n)$ in $\ell^2$. This in turn is linked to a Carleson measure property, as we 
explain below.
By considering the stability of the Carleson measure property under perturbations we are 
able to give a result about exact controllability for perturbed systems. 

Let $p(\lambda,\mu)$ denote the pseudo-hyperbolic metric, given by
\[
p(\lambda,\mu) = \left |\frac{\lambda-\mu}{\lambda+\overline \mu}
\right| \qquad (\lambda, \mu \in \CC_-).
\]
We also require the hyperbolic metric, $d(\lambda,\mu)=2\tanh^{-1} p(\lambda,\mu)$; see, for example
\cite{BM}.

\begin{theorem}\label{thm:perturbdiag}
Suppose that a diagonal system with eigenvalues $(\lambda_n)$ is exactly controllable by
the rank-one controller with parameters $(b_n)$, and that a perturbed system
is also diagonal with eigenvalues $(\mu_n)$ and the same Riesz basis of eigenvectors. Let $p_{n,k}=p(\lambda_n,\lambda_k)$ and $\epsilon_{n}=p(\lambda_n,\mu_n)$.
Provided that (i)~$\sup_n \epsilon_n < 1$, 
 (ii)~$\epsilon_n+\epsilon_k < p_{n,k}$ for all $n \ne k$, 
and (iii)~$\sup_n \sum_{k \ne n} (\epsilon_n+\epsilon_k)|p_{n,k}-1/p_{n,k}| < \infty$, then the perturbed system is also
exactly controllable with the same controller.

If the original system is null-controllable in time $\tau>0$ and conditions (i)--(iii) are satisfied, then the
perturbed system is null-controllable in time $\tau$ with the same controller.
\end{theorem}

\beginpf
By \cite[Thm. 3.1]{JP06},
exact controllability is equivalent to the condition that
\[
\nu_\lambda:=\sum_n  \frac{|\re \lambda_n|^2}{|b_n|^2 \prod_{k \ne n} p(\lambda_n,\lambda_k)^2}\delta_{-\lambda_n}
\]
be a Carleson measure. Here $\delta_{-\lambda_n}$ is a Dirac measure (point mass) at $-\lambda_n$. We show that the analogously defined $\nu_\mu$ is also a Carleson
measure.

Note that if $\sup_n p(\lambda_n,\mu_n)< 1$, then there exists an $\alpha>0$ such that
\[
(1+\alpha)^{-1}|\re \lambda_n| \le |\re \mu_n|\le (1+\alpha) |\re \lambda_n|
\]
 and $|\im \mu_n - \im \lambda_n |\le \alpha |\re \lambda_n|$
for all $n$. (To see this, note that the pseudo-hyperbolic metric is invariant under
translations by $iy$ for $y \in \RR$ and under positive real scalings. Then we may
assume without loss of generality that $\lambda_n=-1$ and thus consider circles
$|(\mu+1)/(\mu-1)|=R < 1$.)

Thus if $|\re \mu_n| \le h$ and $|\im \mu_n-c| \le h/2$, we have $|\re \lambda_n| \le (1+\alpha)h$
and 
\[
|\im \lambda_n - c | \le h/2 + \alpha |\re \lambda_n| \le h/2 + \alpha(1+\alpha)h,
\]
so if $-\mu_n$ is in a fixed Carleson square of size $h$, then $-\lambda_n$ is in a fixed Carleson
square of size at most $(1+2\alpha(1+\alpha))h$.\\

Writing $d_n=d(\lambda_n,\mu_n)=2\tanh^{-1}\epsilon_n$, we have
the inequality
$d(\mu_n,\mu_k) \ge d(\lambda_n,\lambda_k)-(d_n+d_k)$ so 
\begin{eqnarray*}
p(\mu_n,\mu_k) = \tanh \frac12 d(\mu_n,\mu_k) &\ge &
\frac{\tanh \frac12 d(\lambda_n,\lambda_k)-\tanh \frac12 (d_n+d_k)}{1-\tanh \frac12 d(\lambda_n,\lambda_k)\tanh \frac12 (d_n+d_k)} \\
&\ge& p_{n,k}\frac{1-(\epsilon_n+\epsilon_k)/p_{n,k}}{1-p_{n,k}(\epsilon_n+\epsilon_k)},
\end{eqnarray*}
since $\tanh \frac12 (d_n+d_k) \le \tanh \frac12 d_n + \tanh \frac12 d_k$.
Thus convergence of $\prod_{k \ne n} p(\mu_n,\mu_k)$ to a quantity bounded below by a multiple of $\prod_{k \ne n} p(\lambda_n,\lambda_k)$
is assured if 
\[
\sup_n \sum_{k \ne n} (\epsilon_n+\epsilon_k)|p_{n,k}-1/p_{n,k}| < \infty.
\]

The result on null controllability is proved in an identical fashion.

\endpf

\begin{example}
Consider the heat equation, which may be regarded as a diagonal system with eigenvalues $(-n^2)_{n \ge 1}$.
Let us estimate
\[
\sum_{k \ne n} |p_{n,k}-1/p_{n,k}| = \sum_{k \ne n} \frac{4k^2n^2}{|n^4-k^4|}.
\]
We then note that there are absolute constants $C_1, C_2$ and $C_3$ such that
\begin{eqnarray*}
\sum_{k \ne n} \frac{k^2}{|n^4-k^4|} &\le& \left( \sum_{1 \le k < n-\sqrt n}+ \sum_{1 \le |k-n| \le \sqrt n}+ \sum_{k > n+\sqrt n} \right) \frac{k^2}{|n^4-k^4|} \\
& \le & C_1 n \times \frac{n^2}{n^{7/2}} + C_2 \sqrt n \times \frac{n^2}{n^3} + C_3 \sum_{k > n+\sqrt n} \frac{k^2}{k^{7/2}} = O(n^{-1/2}).
\end{eqnarray*}
Thus perturbations asymptotically of size $\epsilon_n = O(n^{-3/2})$, i.e., $|\lambda_n-\mu_n|=O(n^{1/2})$, will allow the conditions of Theorem \ref{thm:perturbdiag} to be satisfied.
That is, if the original system is null-controllable in time $\tau>0$, then so is the perturbed system.
\end{example}

\begin{example}
We may also consider a version of the (damped) wave equation, with eigenvalues
$(-1+in)_{n \in {\mathbb Z}}$. 
A routine calculation shows that for $k \ne n$ we have
\[
|p_{n,k}-1/p_{n,k}| = \frac{4}{|n-k|(4+(n-k)^2)^{1/2}}.
\]
Summing this over $k \ne n$ gives a quantity that is bounded above and below by absolute
constants independent of $n$. Thus provided the supremum of the $\epsilon_n$ is sufficiently
small, the perturbed system inherits the controllability properties of the original system.
\end{example}


\section{Controllability for systems with Jordan blocks}\label{sec:new3} 

So far we have been looking at controllability questions in the case that $A$ is diagonal, but let us now suppose we have a Jordan block $J$ of size $N$ with eigenvalue $\mu$, i.e.,
\[
J=\begin{pmatrix} \mu & 1 & 0 & \ldots & 0 \\
0 & \mu & 1 & \ddots & \vdots \\
\vdots & \ddots & \ddots & \ddots & 0 \\
0 & \ldots & \ddots & \mu & 1 \\
0 & \ldots & \ldots & 0 & \mu
\end{pmatrix},
\]
and thus
\[
 e^{Jt}  
= e^{\mu t}
\begin{pmatrix} 1 & t & t^2/2! & \ldots & t^{N-1}/(N-1)! \\
0 & 1 & t & \ddots & \vdots \\
\vdots & \ddots & \ddots & \ddots & t^2/2! \\
0 & \ldots & \ddots & 1 & t \\
0 & \ldots & \ldots & 0 & 1
\end{pmatrix}.
\]
Therefore, if $b_1,\ldots,b_N$ are the components of $b$ that correspond to the generalized eigenvectors,
then, for exact controllability,
we have to be able to solve a sequence of interpolation conditions of the form
\begin{eqnarray*}
b_1 \hat u(-\mu)-b_2\hat u'(-\mu) + \ldots + \frac{(-1)^{N-1}}{(N-1)!} b_N \hat u^{(N-1)} (-\mu) &=& c_1, \\
\vdots &=& \vdots \\
b_{N-1}\hat u(-\mu) - b_N \hat u'(-\mu) &=& c_{N-1}, \\
b_N \hat u(-\mu)&=& c_N.
\end{eqnarray*}
However, it is clear that if there is more than one Jordan block for any eigenvalue $\mu$, then 
exact controllability (and even null controllability) will be impossible.
\begin{lemma}
If for every $(c_n) \in \ell^2$ there is a function $f \in H^2(\CC_+)$ such that $b_n f(-\lambda_n)=c_n$ for a Blaschke
set $(\lambda_n)$ in $\CC_-$, then there also exists a function $g \in H^2(\CC_+)$ which
satisfies the above interpolation conditions, together with a finite number of additional
interpolation conditions of the form
\[
a_1 g(-\mu_1)=d_1, \quad a_2 g'(-\mu_1)=d_2, \quad \ldots, \quad a_r g^{(r-1)} (-\mu_1) = d_r,
\]
where $a_1,\ldots,a_r$ are nonzero, and $\mu_1$ is not in the set $(\lambda_n)$. 
The same holds if we have additional interpolation conditions of the same form at finitely-many other points not in the set $(\lambda_n)$.
\end{lemma}
\beginpf
For a single point $-\mu_1$, we define
\[
g(s)=f(s)+\frac{p_1B(s)}{s+1}+\frac{p_2 B(s)(s+\mu_1)}{(s+1)^2} + \ldots + \frac{p_r B(s) (s+\mu_1)^{r-1} } {(s+1)^r},
\]
where $p_1,\ldots,p_r$ are constant, and $B$ is the Blaschke product with zeroes $(-\lambda_n)$. We first
choose $p_1$ to obtain the correct value of $g(-\mu_1)$, then $p_2$ to obtain $g'(-\mu_1)$, and so on.
If we   have a further point $\mu_2$ at which we wish to interpolate the function and some of its derivatives, then we add on an
analogous sequence of terms, replacing $B$ by a Blaschke product that has the $(-\lambda_n)$ and $-\mu_1$ as
zeroes (with sufficient multiplicity that the previous interpolation conditions are still satisfied). And so on.
\endpf

The conclusion is that provided the generalized eigenvectors form a Riesz basis, we can ignore
finitely-many repeated eigenvalues (if they have only one Jordan block each), and test the Carleson measure condition for the 
remaining ones, as was done in \cite{JP06}.

Thus we get a more general result than the main result of \cite{JP06}.
\begin{theorem}\label{thm:excon}
A system (\ref{eq:state}) such that $A$ has a Riesz basis of generalized eigenvectors, with all but finitely many of the eigenvalues
$(\lambda_n)$ having algebraic multiplicity 1, and none having a repeated eigenvalue (i.e., only one Jordan
block per eigenvalue, and all but finitely many blocks of size 1) is exactly controllable if and only
if 
\[
\nu:=
\sum_{n=1}^\infty {{|\re \lambda_n|^2}\over{|b_n|^2 d_n^2}}
 \delta_{-\lambda_n}
\]
is a Carleson measure on $\CC_+$, where we take each eigenvalue exactly once, and
\[
d_n=\prod_{k \ne n}\left|\frac{ \lambda_k-\lambda_n}{\lambda_k+ \overline{\lambda_n}}\right| \qquad \hbox{for each} \quad n.
\]
\end{theorem}
Note the general result that if a Riesz basic sequence (in our case, the eigenvectors of $A$) 
spans a subspace $\mathcal M$ of finite codimension, and another
finite sequence is a basis for a (not necessarily orthogonal) complement $\mathcal N$ of $\mathcal M$, then the union of the two sequences is always a
Riesz basis for the whole space. We leave this as an easy exercise.\\

\begin{remark}
We have been discussing exact controllability, but the same techniques can be used to
discuss null-controllability. Under the hypotheses of Theorem
\ref{thm:excon} the necessary and sufficient condition for null controllability in time $\tau>0$ is that
\[
\nu:=
\sum_{n=1}^\infty {{|{\rm Re}\,\lambda_n|^2}\over{|b_n|^2 d_n^2}}e^{2\tau \re \lambda_n}
 \delta_{-\lambda_n}
\]
is a Carleson measure on $\CC_+$. This again extends a similar result for diagonal systems in \cite{JP06}.
\end{remark}

\begin{remark}
In \cite{JPP10}, the theory of interpolation in model spaces was used to study finite-time controllability
for diagonal systems. It would be interesting to adapt these ideas for systems with nontrivial
Jordan blocks.
\end{remark}


\section{Weighted multiple interpolation}\label{sec:new4}

In the previous section, it was shown that additional higher-order interpolation conditions \emph{in a finite number of points} can always be solved for weighted interpolation in $H^2(\CC_+)$,  if the original
interpolation problem can always be solved. In the current section, we want to consider higher-order interpolation conditions  \emph{in infinitely many points}  for weighted interpolation in $H^2(\CC_+)$.
Our results extend the classical results of Vasyunin on multiple interpolation with standard weights, see e.g. \cite{nik}, Lecture X. As in the case of simple weighted interpolation, 
the geometry of the sequence  $(\lambda_n)$ of points which we wish to interpolate, and the sequence of weights  $(G_n)$, are being combined into a Carleson measure
condition, where ``large" weights can compensate for a ``bad" geometry of the sequence  $(\lambda_n)$.   This makes this general weighted setting somewhat more flexible
than the classical setting with standard weights.

More precisely, let $(\lambda_n,K_n)_{n \in \NN} $ be a sequence in $\CC_+ \times \NN$, such that the generalised Blaschke condition
$$
    \sum_{n=1}^\infty K_n \frac{\re \lambda_n}{1+|\lambda_n|^2} < \infty
$$
holds. For each $n \in \NN$, let $G_n$ be an invertible $K_n \times K_n$ matrix.

We want to solve the following weighted multiple interpolation question:
\vspace{0.5cm}
For each sequence $c= (c_n)=  (c_{n,k})_{n \in \NN, k=0, \dots, K_n-1} \in \bigoplus_n \CC^{K_n}$, find a function $f \in H^2(\CC_+)$ such that 
\begin{equation}  \label{eq:freeint}
      G_n \f(\lambda_n)  = c_{n} \qquad n \in \NN
\end{equation}
where $\f(\lambda_n)$ denotes the vector $(f(\lambda_n), f'(\lambda_n), \dots, f^{(K_n-1)}(\lambda_n) ) \in \CC^{K_n}$.

\vspace{0.5cm}

The motivation to consider general invertible $K_n \times K_n$ matrices $G_n$ here stems from the previous section.

Such multiple interpolation questions on $H^2$ have been considered in the work of Vasyunin and others (see e.g. \cite{nik}, Lecture X)
 for the case of the standard weight matrices $G_n = \diag( \frac{(\re\lambda_n)^{k+1/2}}{k!}   )_{k=0, \dots, K_n-1} $. In this case, the evaluation
operator
$$
      H^2(\CC_+)  \rightarrow \bigoplus_n \CC^{K_n}, \qquad f \mapsto (G_n \f(\lambda_n))
$$
is automatically bounded, if it is surjective. However, we do not ask for boundedness of the evaluation operator in our setting, and this allows for much less restrictive conditions.

As usual, for $\lambda \in \CC_+$, let
$$
   b_\lambda (z) = \frac{z- \lambda}{z + \bar \lambda}       \qquad (z \in \CC_+)
$$
denote the Blaschke factor associated to $\lambda$.

Here is the main result of this section:
\begin{theorem}    \label{thm:main}
Let $(\lambda_n,K_n)_{n \in \NN} $ and $(G_n)_{n \in \NN}$ as above, and assume that $(K_n)$ is bounded.
 Then the interpolation problem (\ref{eq:freeint}) has a solution for each  $(c_{n}) \in \bigoplus_n \CC^{K_n}     $,
 if the measure
\begin{equation}    \label{eq:carl}
   \sum_{n=1}^\infty \delta_{\lambda_n} \frac{\re \lambda_n}{\beta_n^{2}} \left  \|  \diag \left ( \frac{(\re\lambda_n)^{k+1/2}}{k!} \right  )_{k=0, \dots, K_n-1}          G_n^{-1} \right  \|^2 
\end{equation}
is a Carleson measure, where 
$$
        \beta_n = \left| \prod_{l=1, l\neq n}^\infty  b_{\lambda_l}^{K_l}(\lambda_n)   \right|.
$$
\end{theorem}
\begin{remark}
In case $G_n =    \diag( \frac{(\re\lambda_n)^{k+1/2}}{k!}   )_{k=0, \dots, K_n-1} $, this is equivalent to the Carleson-Vasyunin
condition $\beta_n \ge \delta$ $(n \in \NN)$ for some
$\delta >0$.
\end{remark}
\proof We use McPhail's method from \cite{mcphail}. By a weak$^*$ compactness argument, it is sufficient to solve the problem for finite sequences $(\lambda_n,K_n)_{1\le n \le N } $, with uniform control 
of norms. So let us solve the problem (\ref{eq:freeint}) for $n=1, \dots, N$.

For $n \le N$, let
$$
  \B_n = \B_n^N =   (\B_n(z))_{j,k} =  \left(     b^k_{\lambda_n}(z)   \prod_{l=1, l\neq n}^N  b_{\lambda_l}^{K_l}(z)
\right)^{(j)} \text{ for } k,j=0, \dots, K_{n}-1,
$$
suppressing the index $N$.
We write
$$
B_{n,k} =  b^k_{\lambda_n}(z)   \prod_{l=1, l\neq n}^\infty  b_{\lambda_l}^{K_l}(z)   \quad \text{ for } k=0, \dots, K_{n} -1,
$$
for the entries in the $0$th row of the matrix $K_n \times K_n$ matrix $\B_n$.

Given $c= (c_n)_{1 \le n \le N} = (c_{n,k})_{1\le n \le N, k=0, \dots, K_n -1} \in  \bigoplus_n^N \CC^{K_n}       $, let
$$
   F_c(z) = \sum_{n=1}^N     e_n(z)  \B_n(z)      \B_n(\lambda_n)^{-1}       G_n^{-1}     c_n    \qquad (z \in \CC_+),
$$
where for each $n$, $e_n \in H^2(\CC_+)$ is a function with $e_n(\lambda_m) =\delta_{nm}$, $e_n^{(j)}(\lambda_m)=0$ for $j=1, \dots, K_n-1$. 

We think of each $\CC^{K_n}$ as a subspace of $l^2$, and of $F_c$ as a function taking values in $l^2$.

Let $f_c$ denote the $0$th component of the vector-valued function $F_c$,
$f_c(z) = (F_c)_0(z)$. Notice that since the $j$th row of $\B_n$ is the $j$th derivative of the $0$th row of $\B_n$ for $0 \le j <K_n$, we have 
\begin{equation}
   F_c(\lambda_n) = (f_c(\lambda_n), f_c'(\lambda_n), \dots,  f_c^{(K_n-1)}(\lambda_n), 0, \dots )
\end{equation}
\begin{lemma}
\begin{equation}   \label{eq:interpol}
  G_n         (f_c(\lambda_n), f_c'(\lambda_n), \dots,  f_c^{(K_n-1)}(\lambda_n))                      =        G_n  F_c(\lambda_n)   = c_n \quad \text{ for }  n=1, \dots, N.
\end{equation}
\end{lemma}
\proof of the lemma. 
This follows directly from the definition of $F_c$.

\qed

Thus 
\begin{equation}
f_c(z) = (F_c(z))_0 =    \sum_{n=1}^N   \sum_{k=0}^{K_n -1}   e_n(z)   B_{n,k}(z)    (    \B_n(\lambda_n)^{-1}     G_n^{-1}     c_{n}   )_k    
\end{equation}
solves the interpolation problem (\ref{eq:freeint}) for $n=1, \dots, N$. We now seek a solution with optimal norm. 

Let $B= B_N = \prod_{l=1}^N  b_{\lambda_l}^{K_l}(z)  $.
Clearly the difference of any two solutions of 
      \begin{equation}   \label{eq:finiteint}
      G_n        (f(\lambda_n), f'(\lambda_n), \dots,  f^{(K_n-1)}(\lambda_n))          = c_{n} \qquad n =1, \dots, N
\end{equation}
in $H^2(\CC_+)$ is contained in $B H^2(\CC_+)$, thus the minimal norm of a solution for a given $c= (c_{n,k}) \in    l^2           $ is
$$
 \inf_{h \in H^2(\CC_+)}      \| f_c(z) - B(z) h(z)   \|_2,   
$$
and we have to estimate
\begin{eqnarray*}
      m_N &=& \sup_{ \|c\|_{l^2} =1}  \inf_{h \in H^2(\CC_+)}      \| f_c(z) - B(z) h(z)   \|_2      \\ 
        &=& \sup_{\|c\|_{l^2} =1}     \| P^- (    \bar B f_c  )\|_2   \\
        &=&  \sup_{h \in H^2(\CC_+), \|h\|_2 =1}    \sup_{ \|c\|_{l^2} =1}    | \langle   \bar B f_c,  \bar h \rangle| . \\
 \end{eqnarray*}

      Write $\tilde c_l =       \B_l(\lambda_l)^{-1}  G_l^{-1} c_{l}$ for $l=1, \dots, N$. We obtain
\begin{eqnarray}         \label{eq:geninter}
 &&  m_N   \nonumber  \\
     &=&  \sup_{h \in H^2(\CC_+), \|h\|_2 =1}    \sup_{ \|c\|_{l^2} =1}  
                 \left|   \sum_{l=1}^N   \sum_{k=0}^{K_l -1}   \left \langle     e_l(z)        \bar b_{\lambda_l}^{K_l-k}(z)  ,  \bar h \right  \rangle        \tilde c_{l,k}                                          \right|  
                   \nonumber \\
         &=&  \sup_{h \in H^2(\CC_+), \|h\|_2 =1}    \sup_{ \|c\|_{l^2} =1}  
                \left|   \sum_{l=1}^N   \sum_{k=0}^{K_l -1}   \left  \langle       \left( \frac{1}{-\bar z - \lambda_l} \right)^{K_l-k}   ,    (\bar z + \lambda_l)^{K_l-k}        
                                                                                                                                        \bar e_l(z)  \bar h \right  \rangle         \tilde c_{l,k}                                 \right|    \nonumber   \\       
          &=&{\sqrt{\pi}}   \sup_{ \|h\|_2 =1}    \sup_{ \|c\|_{l^2} =1}  
                 \left|   \sum_{l=1}^N   \sum_{k=0}^{K_l -1}      \frac{1}{(K_l-k-1)!}      (  ( z + \bar \lambda_l)^{K_l-k}   h)^{(K_l-k-1)}(\lambda_l)   \tilde c_{l,k}              \right| \nonumber  \\         
              &=& {\sqrt{\pi}} \sup_{\|h\|_2 =1}    \sup_{ \|c\|_{l^2} =1}  
   \left|   \sum_{l=1}^N   \left \langle   \left( \frac{1}{(K_l-k-1)!}      (  ( z + \bar \lambda_l)^{K_l-k}   h)^{(K_l-k-1)}(\lambda_l)   \right)_{k=0, \dots, K_n-1} ,
         \overline{ \B_l(\lambda_l)^{-1}     G_{l}^{-1}          c_{l}  } \right \rangle   \right|   \nonumber \\
                      &=& {\sqrt{\pi}}  \sup_{ \|h\|_2 =1}   
                 \left(   \sum_{l=1}^N   \|     (G_{l}^{-1} )^t       (\B_l(\lambda_l)^{-1})^t         \left( \frac{1}{(K_l-k-1)!}      (  ( z + \bar \lambda_l)^{K_l-k}   h)^{(K_l-k-1)}(\lambda_l)   \right)_{k=0, \dots, K_n-1}     \|_2^2  \right)^{1/2}     \nonumber \\        
                  &=& {\sqrt{\pi}}  \sup_{h \in H^2(\CC_+), \|h\|_2 =1}   
                 \left(   \sum_{l=1}^N  \frac{1}{\re \lambda_l}  \|   (  G_l^{-1})^t  (\B_l(\lambda_l)^{-1})^t   \left(     \langle   h (z + \bar \lambda_l),  e_{\lambda_l,     K_l-k-1} \rangle     \right)_{k=0, \dots, K_n-1}     \|_2^2  \right)^{1/2} 
                            \end{eqnarray}
  where the supremum in the last line is taken on the dense subspace of functions with $ h (z+1) \in H^2$.                         
 Here, we have used the Malmquist--Walsh functionals
 $$
          e_{\lambda,j} =   \frac{1}{\sqrt{\pi}} b_\lambda(z)^{j}   \frac{(\re \lambda)^{1/2}}{z + \bar \lambda}   \qquad (j \ge 0, \lambda \in \CC_+)
 $$
 and the elementary identity
 $$
        \langle u, (-1)^j \frac{1}{\sqrt{\pi}(z + \bar \lambda)^{j+1}} \rangle =\frac{1}{j!} u^{(j)}(\lambda)   \text{ for } u \in H^2(\CC_+), \lambda \in \CC_+, j\in \NN_0.
 $$

For  each $\lambda \in \CC_+$, $K \in \NN_0$, the Malmquist--Walsh functionals $(e_{\lambda,j})_{j=0, \dots, K-1}$ form an orthonormal basis of the model space $K_{b_\lambda^K} = H^2 \ominus b_{\lambda}^K H^2$
(see e.g. \cite{nik}). 
We write
$$
P_{\lambda, K} :H^2 \rightarrow H^2 \ominus b_\lambda^K H^2 \simeq \CC^K, \qquad P_{\lambda,K} h =   
           \sum_{j=0}^{K-1}    \langle h, e_{\lambda_l, K-j -1} \rangle e_{\lambda_l, K-j -1}  
$$
for the orthogonal projection onto the model space in the basis of the Malmquist--Walsh functionals. 

 Let us remark here that the identity (\ref{eq:geninter})
proven above  characterises the minimal norm solutions of the interpolation problem (\ref{eq:finiteint}), and does not depend on the bound of the sequence $(K_n)$:

\begin{theorem}   \label{thm:char}
The interpolation problem (\ref{eq:freeint}) has a solution for each $(c_{n,k}) \in l^2$, if and only if the operator
\begin{eqnarray*}
   \J: && H^2 \rightarrow \bigoplus_n  \CC^{K_n}, \\
      (\J h)_n &= & (G_n^{-1})^t   (\B_n(\lambda_n)^{-1})^t  \left( \frac{1}{(K_n-k-1)!}      (  ( z + \bar \lambda_n)^{K_n-k}   h)^{(K_n-k-1)}(\lambda_n)   \right)_{k=0, \dots, K_n-1}   \\
     &= & \frac{1}{(\re \lambda_n)^{1/2}} (G_n^{-1})^t   (\B_n(\lambda_n)^{-1})^t  \langle    (  ( z + \bar \lambda_n) h, e_{\lambda_n, K_n-k -1}  \rangle)_{k=0, \dots, K_n-1}.   \\
      &= &\frac{1}{(\re \lambda_n)^{1/2}} (G_n^{-1})^t   (\B_n(\lambda_n)^{-1})^t    P_{\lambda_n, K_n}        (  ( z + \bar \lambda_n) h)  \\
\end{eqnarray*}
is bounded. In this case, the norm of the operator $\J$ coincides with the minimal interpolation norm,
$$
   m=   \sup_{ \|c\|_{l^2} =1}  \inf \{    \| f  \|_2    :     f \in H^2(\CC_+), G_n \mathbf{f} (\lambda_n) = c_n \text{ for all } n \in \NN\} = \| \J\|.
$$
\end{theorem}
\proof of Theorem \ref{thm:char}.
This has already been proved for finitely many points and a  weak$^*$ compactness argument gives the result in general.

\qed

%
                 
 We return to the proof of Theorem \ref{thm:main}, controlling the operator $\J$ in Theorem \ref{thm:char}  by a Carleson embedding.
 This requires boundedness of the sequence $(K_n)$. From now on, we write $\lesssim$, $\gtrsim$ for estimates with constants depending only on $\sup_{n}K_n$.
 
 Cauchy's integral inequalities easily imply
 \begin{multline*}
      |  (  ( z + \bar \lambda_l)^{K_l-k}   h)^{(K_l-k-1)}(\lambda_l)| \lesssim      \frac{ (K_l-k-1)! }{(\re \lambda_l)^{K_l-k-1}} (\re \lambda_l)^{K_l-k}   \sup_{|z- \lambda_l| \le \re \lambda_l/4} |h(z)|  \\
      = 
          (K_l-k-1)!  \re \lambda_l \sup_{|z- \lambda_l| \le  \re \lambda_l/4} |h(z)|.
 \end{multline*}      
 
 If $ \sum_{n=1}^\infty \delta_{\lambda_n} (\re \lambda_n)^2 \| \B_n(\lambda_n)^{-1}G_n^{-1} \|^2  $ is a Carleson measure, then 
 $$
        \sum_{n=1}^\infty \delta_{\mu_n} (\re \lambda_n)^2 \|\B_n(\lambda_n)^{-1}G_n^{-1} \|^2    
 $$     
 is also a Carleson measure for any sequence $(\mu_n)$ in $\CC_+$ with $|\mu_n - \lambda_n | \le \re \lambda_n/4   $ for all $n$. In this case, 
 its Carleson constant satisfies
  $$
       \carl(\sum_{n=1}^\infty \delta_{\mu_n} (\re \mu_n)^2 \| \B_n(\lambda_n)^{-1} G_n^{-1}\|^2)   \le 
              2  \carl(\sum_{n=1}^\infty \delta_{\lambda_n} (\re \lambda_n)^2 \|\B_n(\lambda_n)^{-1}    G_n^{-1}\|^2   ).
 $$
 Thus by (\ref{eq:geninter}),
 \begin{multline*}
        m_N \lesssim  \sup_{(\mu_n),  |\mu_n  - \lambda_n| < \re \lambda_n/4 }   \,    \sup_{h \in H^2(\CC_+), \|h\|_2 =1}    
                 \left(   \sum_{l=1}^N   \|  (G_n^{-1})^t     (\B_l(\lambda_l)^{-1})^t   \|^2     (\re \lambda_n)^2     |h(\mu_n)|^2     \right)^{1/2} \\
                   \lesssim
                      \carl(\sum_{n=1}^\infty \delta_{\lambda_n} ( \re \lambda_n)^2  \|\B_n(\lambda_n)^{-1}   G_n(\lambda_n)^{-1} \|^2 )^{1/2}.
 \end{multline*}
 
 To finish the proof of the sufficiency of the Carleson measure condition, we only need to remark that the matrices 
 \begin{equation}  \label{eq:matrix}
    \frac{1}{ \beta_n}         \diag( \frac{(2\re\lambda_n)^{k}}{k!}   )_{k=0, \dots, K_n-1}          \B_n(\lambda_n) 
 \end{equation}
 are bounded above and below, with bounds only depending on $K_n$. To see this, note that the matrix  (\ref{eq:matrix})   
 is superdiagonal with all diagonal entries being equal to $1$, and that the remaining entries are bounded with a bound only depending on $K_n$.
 
 \vspace{0.5cm}

   \qed

  {\bf Remark.} In the case of multiple interpolation with standard weights, it is not difficult to see by testing on reproducing kernels that the Carleson condition in Theorem \ref{thm:main} is also
  necessary. In the general case, it appears that solvability of the weighted multiple interpolation problem (\ref{eq:freeint}), 
  or equivalently, the boundedness of the embedding $\J$ in Theorem \ref{thm:char}, cannot be characterised by a simple
  Carleson condition.

\section{Admissibility for contraction semigroups}\label{sec:new5}

In this section we exploit the positive solution to the Weiss conjecture for 
contraction semigroups and finite-rank observation operators \cite{JP01}
in order to deduce a result on the admissibility of observation
operators for perturbed systems. Because of the well-known duality between 
admissibility for control
and observation operators (see, e.g., \cite{JP04}) one can deduce a similar result
for admissibility of control operators, which we shall not state separately.

 Recall that if $(T(t))_{t \ge 0}$ is a $C_0$ semigroup
on a Hilbert space with infinitesimal generator $A$, then an $A$-bounded operator
$C: D(A) \to Y$ is said to be {\em (infinite-time) admissible\/} if there is a constant $K>0$ such that
\[
\int_0^\infty \|CT(t)x\|^2 \, dt \le K \|x\|^2 \qquad (x \in D(A)).
\]
For more on this concept we refer to the survey \cite{JP04} and the book \cite{TW}.\\

A key perturbation result can be found in \cite[Thm. 5.4.2]{TW}, which is based on
\cite{weiss94,HW}. The idea here is to consider a closed-loop semigroup.

\begin{theorem}[Hansen and Weiss]
Suppose that $B \in \LL(Z,H)$ and $C:D(A) \to Z$ is an admissible
observation operator for $(T(t))_{t \ge 0}$. Then the operator
$A+BC$ generates a $C_0$ semigroup $(T^{cl}(t))_{t \ge 0}$, such that, for every
Hilbert space $Y$ the space of admissible observation operators for $(T(t))_{t \ge 0}$ mapping into $Y$ is the same
as the corresponding space for $(T^{cl}(t))_{t \ge 0}$.
\end{theorem}

Our approach is rather different, and does not involve admissible perturbations.

Recall that $A$-boundedness means that there exist constants $a,b > 0$ such that
\beq\label{eq:abdd}
\|Cx\| \le a\|Ax\|+b \|x\| \qquad (x \in D(A)),
\eeq
and we write $a_0$ for the infimum of all $a \ge 0$ such that
(\ref{eq:abdd}) holds for some $b >0$. See for example \cite[p.~169]{EN}.

The following lemma links $A$-boundedness with resolvent conditions.
\begin{lemma}\label{lem:resolvent}
Let $\Delta: D(A) \to Y$ satisfy $\|\Delta(s-A)^{-1}\|=R$ for some $s \in \CC$. Then $\Delta$ is $A$-bounded with
constant $a_0 \le R$.

Conversely, if $A$ generates an analytic semigroup and $\Delta$ has $A$-bound $0$, then for all $\varepsilon>0$
there exists $w_0>0$ such that
\beq\label{eq:w0bound}
\sup_{\re s>w_0} \|\Delta(s-A)^{-1}\| < \epsilon. 
\eeq
\end{lemma}
\beginpf
For $x \in D(A)$, write $y=(s-A)x$. Then
\begin{eqnarray*}
\|\Delta x\| & = & \|\Delta(s-A)^{-1}y\| \\
& < & R \|y\| = R \|(s-A)x\| \\
& \le & R|s| \|x\|+ R \|Ax\|,
\end{eqnarray*}
establishing the first result. 

Now, if $A$ generates an analytic semigroup, then  there
exists $w_0 > 0$ and $m>0$ such that
\[
\|(s-A)^{-1}\| \le m/|s| \qquad (\re s > w_0),
\]
and so for each $\eta>0$ there exists $b>0$ with
\begin{eqnarray*}
\|\Delta(s-A)^{-1}x\| & \le & \eta \|A(s-A)^{-1}x\|+ b\|(s-A)^{-1}x\| \\
& \le & \eta\|-x+s(s-A)^{-1}x\| + bm \|x\|/|s| \\
& \le & \eta( \|x\|+ m\|x\|) + bm \|x\|/|s|.
\end{eqnarray*}
 Now, letting $\eta \to 0$ and then choosing $s$ with real part sufficiently large, we have
 (\ref{eq:w0bound}).
 \endpf
 
\begin{theorem}
Let $H$ be a Hilbert space and $A$ the generator of a contraction semigroup
$(T(t))_{t \ge 0}$
on $H$, with domain $D(A)$. Suppose that $\Delta: D(A) \to H$ is dissipative and
there exists $w_0>0$ such that
\beq\label{eq:perturb}
\sup_{\re s>w_0} \|\Delta(s-A)^{-1}\| < 1. 
\eeq
Let $Y$ be a finite-dimensional Hilbert space. Then the space of all admissible
observation operators $C:D(A) \to Y$ for $(T(t))_{t \ge 0}$  is equal to the
corresponding space for the semigroup $(T_\Delta(t))_{t \ge 0}$ generated by $A+\Delta$.
\end{theorem}

\beginpf
By Lemma \ref{lem:resolvent}, $\Delta$ is $A$-bounded with $A$-bound $a_0<1$. Now
by \cite[p.~173, Thm.~2.7]{EN}, $A+\Delta: D(A) \to H$ generates a contraction semigroup.

Next, by \cite[Cor.~1.4]{JP01} (see also \cite{JP04}), $C$ is admissible for $(T(t))_{t \ge 0}$ if and only if
there exists $M>0$ and $w \in \RR$ such that
\beq\label{eq:resolvent}
\|C(s-A)^{-1}\| \le \frac{M}{\sqrt{\re s - w}} \qquad \hbox{whenever} \quad \re s > w.
\eeq
Clearly, a similar equivalence holds for $(T_\Delta(t))_{t \ge 0}$, on replacing $A$ by $A+\Delta$.\\

Suppose that $C$ is admissible for $(T(t))_{t \ge 0}$ and that (\ref{eq:resolvent}) holds.
Now 
\[
C(s-A)^{-1}=C(s-A-\Delta)^{-1}(I-\Delta(s-A)^{-1}),
\]
 and hence,  
using (\ref{eq:perturb}), we see that, for $\re s > \max\{w,w_0\}$, we have
\begin{eqnarray*}
\|C(s-A-\Delta)^{-1}\| &\le& \|C(s-A)^{-1}\| \, \|(I-\Delta(s-A)^{-1})^{-1}\| \\
& \le & \frac{M M'}{\sqrt{\re s - w}}
\end{eqnarray*}
for some constant $M'>0$. 
It follows that $C$ is admissible for $(T_\Delta(t))_{t \ge 0}$.\\

For the converse implication, we may use the simpler inequality
\[
\|C(s-A)^{-1}\| \le \|C(s-A-\Delta)^{-1}\| \, \|I + \Delta(s-A)^{-1}\|
\]
to deduce that $C$ is admissible for $(T(t))_{t \ge 0}$ whenever it is admissible for
the perturbed semigroup: indeed, this holds as soon as  $\Delta(s-A)^{-1}$ is bounded, and not 
necessarily a strict contraction.

\endpf

\begin{remark}
The same result holds if $\Delta=\Delta_1+E$, where $\Delta_1$ is a  dissipative $A$-bounded operator with $a_0<1$
and $E$ is an arbitrary bounded operator. For now, using \cite[p.~158, Thm.~1.3]{EN} we see that
$A+\Delta$ generates a semigroup $(T_{\Delta}(t))_{t \ge 0}$ with
growth bound
\[
\|T_\Delta(t)\| \le e^{\|E\|t}.
\]
\end{remark}

\section*{Acknowledgement} This work was supported by EPSRC grant EP/I01621X/1.

\end{document}